\documentclass[11pt, a4paper, reqno]{amsart}
\usepackage[a4paper, margin=3.75cm]{geometry}
\usepackage[T2A]{fontenc}
\usepackage[cp866]{inputenc}
\usepackage[russian,english]{babel}
\usepackage{amsfonts,amssymb,mathrsfs,amscd,comment}
\usepackage{amsthm,graphicx} 
\usepackage[width=1.0\textwidth]{caption}
\def\NoBlackBoxes{\overfullrule0pt}

\NoBlackBoxes
\theoremstyle{plain}
\newtheorem{theorem}{Theorem}
\newtheorem{conjecture}{Conjecture}
\newtheorem{corollary}{Corollary}

\theoremstyle{definition}

\usepackage[breaklinks]{hyperref}
\hypersetup{
unicode=true,
colorlinks=true,
linkcolor=blue,
citecolor=blue,
urlcolor=blue,
filecolor=blue,
bookmarksnumbered=true,
pdfstartview=FitH,
pdfhighlight=/N
}
\theoremstyle{main}
\let\savedef=\endproof
\def\endproof{~$\square$\savedef}
\def\bad{\spaceskip=0.33emplus0.6emminus0.15em\immediate\write5{\string\bad}}

\theoremstyle{plain}

\let\myh\widehat
\let\myo\overline

\def\HH{\mathscr H}

\def\RR{\mathbb R}
\def\CC{\mathbb C}

\def\NN{\mathbb N}

\def\mdeg{\operatorname{deg}}

\def\mcap{\operatorname{cap}}

\def\supp{\operatorname{supp}}
\def\const{\operatorname{const}}
\def\<{\left\langle}
\def\>{\right\rangle}
\def\({\left(}
\def\){\right)}
\def\[{\left[}
\def\]{\right]}
\NoBlackBoxes

\def\bad{\spaceskip=0.33emplus0.6emminus0.15em\immediate\write5{\string\bad}}
\def\NN{\mathbb N}
\def\RR{\mathbb R}
\def\CC{\mathbb C}

\def\zz{\mathbf z}
\def\RS{\mathfrak R}

\def\ff{w}
\def\myint{\operatorname{int}}
\def\sM{\mathscr M}

\let\myh\widehat

\let\myo\overline
\def\const{\operatorname{const}}

\def\({\left(}
\def\){\right)}
\def\[{\left[}
\def\]{\right]}
\def\<{\left\langle}
\def\>{\right\rangle}

\let\leq\leqslant
\begin{document}

\selectlanguage{english}

\title{Maximum Principle and Asymptotic Properties of Hermite--Pad\'e Polynomials}

\author[Sergey~P.~Suetin]{Sergey~P.~Suetin}
\address{Steklov Mathematical Institute of the Russian Academy of Sciences, Moscow, Russia}
\email{suetin@mi-ras.ru}

\date{September 21, 2021}

\maketitle

\markright{Maximum Principle and Hermite--Pad\'e Polynomials}

\begin{abstract}
In the paper, we discuss how it would be possible to succeed in Stahl's novel approach, 1987--1988, to explore Hermite--Pad\'e polynomials based on Riemann surface properties.

In particular, we explore the limit zero distribution of type I Hermite--Pad\'e polynomials $Q_{n,0},Q_{n,1},Q_{n,2}$, $\mdeg{Q_{n,j}}\leq{n}$, for a collection of three analytic elements $[1,f_\infty,f^2_\infty]$. The  element $f_\infty$ is an element of a function $f$ from the class  $\mathbb C(z,\ff)$ where $\ff$ is supposed to be from  the class $Z_{\pm1/2}([-1,1])$ of multivalued analytic functions  generated by the inverse Zhukovskii function with the exponents from the set $\{\pm1/2\}$. The Riemann surface corresponding to $f\in\mathbb C(z,\ff)$ is a four-sheeted Riemann surface  $\RS_4(\ff)$ and all branch points of $f$ are of the first order (i.e., all branch points are of square root type).

It is proved that there exits  the similar limit zero distribution for all three polynomials $Q_{n,j}$. The answer is done in terms of Nuttall's condenser which was introduced by E.~Rakhmanov and the author in 2013. The corresponding limit measure is supported on the second plate of  Nuttall's condenser which coincides with the projection of the boundary between the first and the second sheets of the three-sheeted Riemann surface $\RS_3(\ff_\infty)$ associated by Nuttall with the element $\ff_\infty$. The limit measure is a unique solution of a mixed Green-logarithmic theoretical-potential equilibrium problem.
In general the surface $\RS_3(\ff_\infty)$ can be of any genus.

Since the algebraic function $f\in\mathbb C(z,\ff)$ is of fourth order and we consider the triple of the analytic elements $[1,f_\infty,f^2_\infty]$ but not the quadruple $[1,f_\infty,f^2_\infty,f^3_\infty]$ ones, the result is new and does not follow from the known results.

As in previous  paper \href{https://arxiv.org/abs/2108.00339}{arXiv: 2108.00339} and following to Stahl's ideas, 1987--1988, we do not use the orthogonality relations at all. The proof is based on the maximum principle only.

Bibliography:~\cite{Sue21d}~titles.
\end{abstract}

\setcounter{tocdepth}{1}\tableofcontents

\footnotetext[0]{This work was performed at the Steklov International Mathematical Center and supported by the Ministry of Science and Higher Education of the Russian Federation (agreement no. 075-15-2019-1614).}

\section{Introduction}\label{s1}

\subsection{}\label{s1s1}
Let  $\varphi(z)=z+(z^2-1)^{1/2}$ be the inverse Zhukovskii function; here and in what follows such branch of the function $(z^2-1)^{1/2}$ is taken that outside of the segment $E:=[-1,1]$ it holds $(z^2-1)^{1/2}/z\to1$ as $z\to\infty$. Therefore, we have that $\varphi(z)/z\to2$ as $z\to\infty$.

Let $m\in\NN$ and $A_j,B_j\in\RR$ be real numbers with the following properties: $A_1<B_1<\dots<A_k<B_k<-1$, $1<A_{k+1}<B_{k+1}<\dots<A_m<B_m$. Set
\begin{equation}
\ff(z):=\prod_{j=1}^m\(\frac{A_j-1/\varphi(z)}{B_j-1/\varphi(z)}\)^{1/2}, \quad z\in D:=\myh{\CC}\setminus{E}.
\label{1}
\end{equation}
We denote the class of all analytic functions given by the explicit representation~\eqref{1} by $\mathscr Z_{\pm1/2}(E)$. We emphasize
once again that it is supposed that the parameters $A_j$ and $B_j$ satisfy the above conditions. Notice that for $m=1$ the systems $\ff,\ff^2$ and $\ff,\ff^2,\ff^3$ form Nikishin systems (see~\cite{Sue18c}).

Being defined in such a way, the function $\ff$ is an algebraic function of fourth order. All of its branch points are of the first order (i.e., all branch points are of square root type) and they constitute the set $\Sigma=\Sigma_{\ff}=\{\pm1,a_1,b_1,\dots,a_m,b_m\}$  where $a_j=(A_j+1/A_j)/2$, $b_j=(B_j+1/B_j)/2$,  $j=1,\dots,m$.
The Riemann surface (RS) corresponding to $\ff$ is a four-sheeted RS $\RS_4(\ff)$.
Under the above condition on $\varphi(z)$,  there is an analytic element $\ff_\infty\in\HH(\infty)$ of $\ff$  with the following property: $\ff_\infty(\infty)=\prod\limits_{j=1}^m\sqrt{A_j/B_j}>0$. The element $\ff_\infty$ extends
as a holomorphic (single-valued analytic) function to the domain $D=\myh{\CC}\setminus{E}$. Set
 $F:=\bigsqcup\limits_{j=1}^m[a_j,b_j]$.

Let $f\in\CC(z,\ff)$. Then $f $ is a (single-valued) meromorphic function on the RS $\RS_4(\ff)$.
Therefore the asymptotic properties of type I Hermite--Pad\'e polynomials (HP-polynomials) for the collection of four functions $[1,f,f^2,f^3]$ follows directly from the results of the paper~\cite{KoPaSuCh17} (see also~\cite{Kom21} and~\cite{Kom21b}). But this is not the case for the collection of three functions $[1,f,f^2]$ when $f\in\CC(z,\ff)$. Notice also that in general function $f\in\CC(z,\ff)$ is a complex-valued function on the real line. Therefore a very general and powerful approach initiated by Gonchar and Rakhmanov~\cite{GoRa81} in 1981 and  substantially generalized and advanced in the excellent papers~\cite{GoRaSo97} and ~\cite{ApLy10}, is not applicable in the case under consideration in the current paper.

\subsection{}\label{s1s2}
Let $f\in\CC(z,\ff)$ and $f_\infty\in\HH(\infty)$ be the element of $f$ which corresponds to the element $\ff_\infty$ being described before; in particular $\ff_\infty(\infty)=\prod\limits_{j=1}^m\sqrt{A_j/B_j}>0$. We suppose that this condition holds up to the end of the paper.

Let $n\in\NN$ be fixed and polynomials $Q_{n,0},Q_{n,1},Q_{n,2}$ of degree $\leq{n}$ be defined (not uniquely) from the relation
\begin{equation}
R_n(z):=(Q_{n,0}+Q_{n,1}f_\infty+Q_{n,2}f_\infty^2)(z)=O\(\frac1{z^{2n+2}}\),
\quad z\to\infty.
\label{2}
\end{equation}
Those polynomials $Q_{n,0},Q_{n,1},Q_{n,2}$ are type I Hermite--Pad\'e polynomials for the multiindex $(n,n,n)$ and the collection $[1,f_\infty,f_\infty^2]$. Function $R_n(z)=R_n(z;f)$ is the error function.

For each (positive Borel) measure $\mu$, $\supp{\mu}\subset\CC$, we demote by $V^\mu(z)$ the logarithmic potential of $\mu$:
$$
V^\mu(z):=\int\log\frac1{|z-\zeta|}\,d\mu(\zeta),\quad z\in\CC\setminus\supp{\mu}.
$$
Let $g_F(z,\zeta)$, $z,\zeta\in\Omega:=\myh{\CC}\setminus{F} $ be the Green function for the domain $\Omega$ with logarithmic singularity at the point $z=\zeta$. Let
\begin{equation}
G_F^\mu(z):=\int g_F(z,\zeta)\,d\mu(\zeta)
\label{3}
\end{equation}
be the corresponding Green potential of the measure $\mu$. In a similar way we define the functions $g_E(z,\zeta)$ and $G_E^\mu(z)$ for the domain $D=\myh{\CC}\setminus{E}$.

It is well-known (see~\cite{RaSu13} and the bibliography therein) that there exists a unique probability measure $\lambda_E$, supported on the compact set $E$, $\lambda_E\in M_1(E)$, with the following property\footnote{We suppose that there will be no ambiguity in the notations $w_E$, $w_\infty$ and etc.}:
\begin{equation}
3V^{\lambda_E}(x)+G^{\lambda_E}_F(x)\equiv w_E=\const,
\quad x\in E.
\label{4}
\end{equation}
Let $\lambda_F=\beta_F(\lambda_E)\in M_1(F)$ be the balayage of the measure $\lambda_E$ from the domain $\Omega$ to the compact set  $F=\partial\Omega$.

\begin{figure}
\centerline{
\boxed{
\includegraphics[width=\textwidth]{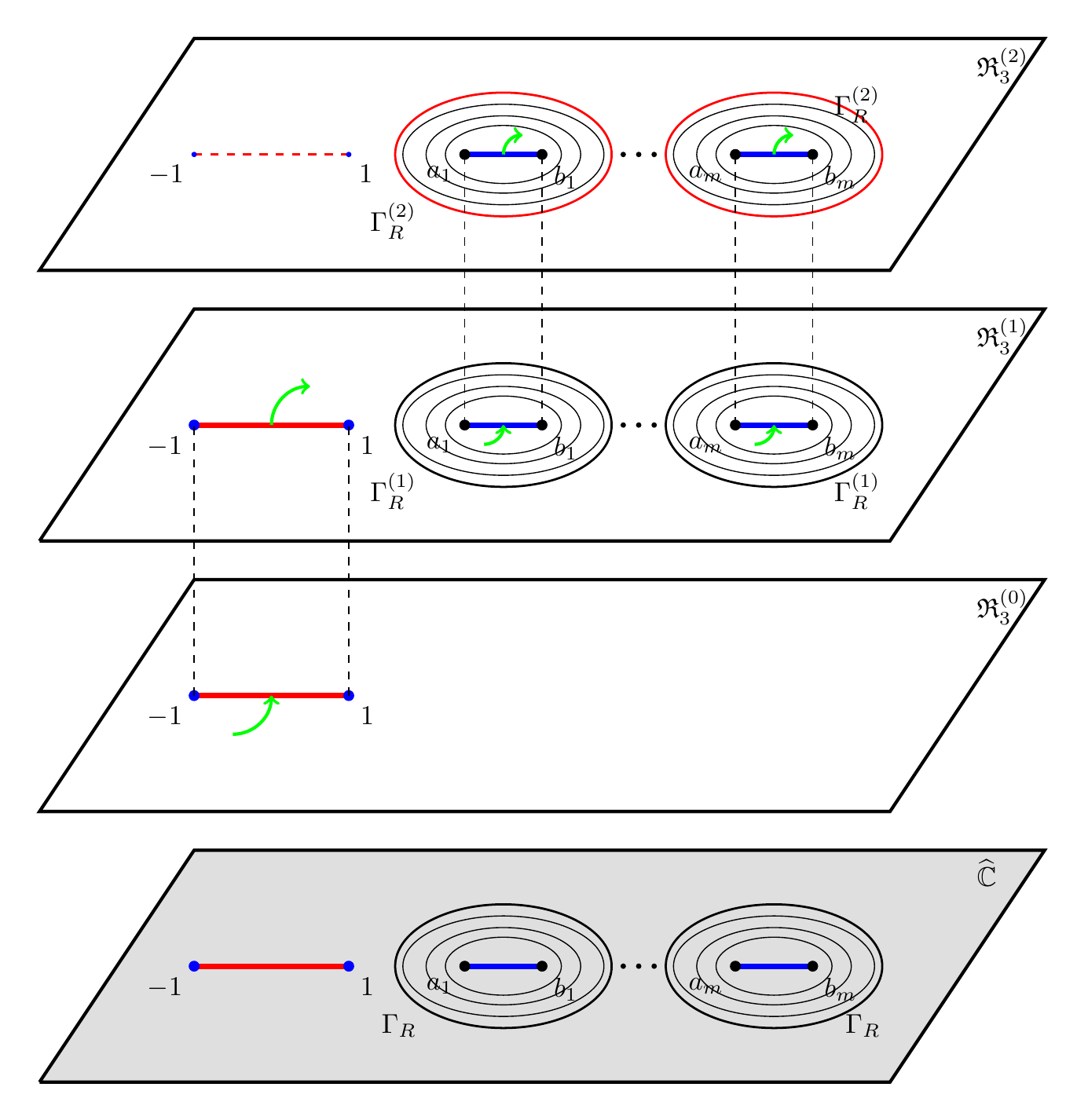}
}}
\caption{On the figure it is represented the three-sheeted Riemann  surface $\RS_3(\ff_\infty)$ associated by Nuttall with the given element $\ff_{\infty^{(0)}}=\ff_\infty\in\HH(\infty)$. The double sided red segment is the boundary between the zero $\RS_3^{(0)}$ and the first $\RS_3^{(1)}$ sheets of RS $\RS_3(\ff_\infty)$. The double-sheeted black segments are the boundary between the first $\RS_3^{(1)}$ and the second $\RS_3^{(2)}$ sheets of RS $\RS_3(\ff_\infty)$. Oval closed curves are the lifting of the level curves of the function $G_{F}^{\lambda_E}(z)$ to the first sheet and to the second sheet from the ``physical'' plane~$\myh{\CC}$. The physical plane itself is colored in grey.
In the case under consideration both plates $E$ and $F$ of the Nuttall condenser are real
compact sets.}
\label{fig_1}
\end{figure}

On the base of the properties of potentials of the measure $\lambda_E$
and according to the scheme proposed in~\cite{RaSu13} (see also~\cite{Sue19}), it is possible to construct a three-sheeted Riemann surface   $\RS_3(f_\infty)$ with the so-called Nuttall partition into open  sheets  $\RS_3^{(0)}\ni z^{(0)},\RS_3^{(1)}\ni z^{(1)}$ and $\RS_3^{(2)}\ni z^{(2)}$, $\pi_3(z^{(j)})=z$, and with the  following properties. The given element $f_\infty\in\HH(\infty)$ is lifted to the infinity point $\infty^{(0)}\in\RS^{(0)}$ and after then extends as a meromorphic function up to the boundary between the first sheet $\RS_3^{(1)}$ and the second sheet $\RS_3^{(2)}$. It is easy to see that this element extends even beyond this boundary a little on the second sheet. The following function
\begin{equation}
\begin{aligned}
u(z^{(2)}):&=G_F^{\lambda_E}(z)-V^{\lambda_E}(z),
\quad z\notin{F},\\
u(z^{(1)}):&=-G_F^{\lambda_E}(z)-V^{\lambda_E}(z),
\quad z\notin E\cup F,\\
u(z^{(0)}):&=2V^{\lambda_E}(z)-w_E,\quad z\notin E,
\end{aligned}
\label{5}
\end{equation}
extends to the whole set $\RS_3(f_\infty)\setminus\{\infty^{(0)},\infty^{(1)},\infty^{(2)}\}$ as a harmonic function $u(\zz)$ with the following singularities at the infinity points:
\begin{equation}
\begin{aligned}
u(z^{(2)}):&=\log|z|+O(1),\quad z\to\infty,\\
u(z^{(1)}):&=\log|z|+O(1),\quad z\to\infty,\\
u(z^{(0)}):&=-2\log|z|+O(1),\quad z\to\infty.
\end{aligned}
\label{6}
\end{equation}
From~\eqref{5} it follows (see~\cite{RaSu13} and subsection~\ref{s2s1} below) that
\begin{equation}
u(z^{(0)})<u(z^{(1)})<u(z^{(2)}).
\label{7}
\end{equation}
Thus from~\eqref{6} and~\eqref{7} it follows that the partition into   three (open) sheets $\RS_3^{(0)},\RS_3^{(1)}$ and $\RS_3^{(2)}$ is indeed the  Nuttall partition of $\RS_3(f_\infty)$ with respect to the point $\infty^{(0)}$. For  more details about the construction of the Riemann surface $\RS_3(f_\infty)$ see~\cite{RaSu13},~\cite{Sue19}, Section~\ref{s3} below and Fig.~\ref{fig_1}. Notice that the structure of the Nuttall partition with respect to the point $\infty^{(0)}$
for the four-sheeted RS $\RS_4(\ff)$ of the function $\ff$  was analyzed in detail in~\cite{IkSu21b}. From that paper it follows that under the above conditions on the parameters $A_j$ and $B_j$ in both cases of $\RS_3(w_\infty)$ and $\RS_4(w_\infty)$ zero and first sheets of the partitions are just the same. But it is not the case in the general situation; for more details see~\cite[\S~4]{IkSu21b} and Section~\ref{s3} below.

For a polynomial $Q\in\CC[z]\setminus0$ let
$$
\chi(Q)=\sum_{\zeta:Q(\zeta)=0}\delta_\zeta
$$
be the corresponding zero counting measure.

The main result of the paper is the following statement (cf.~\cite{Sue18d},~\cite{KoPaSuCh17},~\cite{Kom21},~\cite{Kom21b}).

\begin{theorem}\label{the1}
Let $f\in\CC(z,\ff)$ and the analytic element $f_\infty\in\HH(\infty)$ satisfies the above conditions. Then for type I Hermite--Pad\'e polynomials $Q_{n,j}$, $j=0,1,2$, we have that as $n\to\infty$:
\begin{align}
\frac1n\chi(Q_{n,j})&\overset{*}\longrightarrow\lambda_F,\quad j=0,1,2,
\label{8}\\
\biggl|\frac{Q_{n,1}(z)}{Q_{n,2}(z)}+\bigl(f(z^{(0)})+f(z^{(1)})\bigr)
\biggr|^{1/n}&\overset{\mcap}\longrightarrow e^{-2G_F^{\lambda_E}(z)}<1,\quad z\in\Omega.
\label{9}
\end{align}
\end{theorem}

From the proof of Theorem~\ref{the1} it follows  that
\begin{corollary}\label{cor1}
The sequences $\{n-\mdeg{Q_{n,j}}\}$, $j=0,1,2$, are bounded.
\end{corollary}

Finally we notice that the asymptotic properties of Hermite--Pad\'e polynomials are still of unabated interest (see~\cite{Apt08},~\cite{ApLoMa17},~\cite{LoMeSz19},~\cite{Sor20},~\cite{Lys20},~\cite{Sor20},~\cite{ApLy21} and the bibliography therein). An accessible presentation of the potential theory on Riemann surface is given in the papers~\cite{Chi18}--\cite{Chi20}.

\section{Proof of Theorem~\ref{the1}}\label{s2}

\subsection{}\label{s2s1}
For two positive sequences $\{\alpha_n\}$ and $\{\beta_n\}$ the relation $\alpha_n\asymp \beta_n$ means that  $0<C_1\leq \alpha_n/\beta_n\leq C_2<\infty$ for $n=1,2,\dots$ and some constants $C_1,C_2$ which do not depend on $n$. For two sequences $\{\alpha_n(z)\}$ and $\{\beta_n(z)\}$ of functions holomorphic in a domain $G$ the relation $\alpha_n\asymp\beta_n$ means that for each compact set $K\subset G$ and $n=1,2,\dots$ \  $0<C_1\leq |\alpha_n(z)/\beta_n(z)|\leq C_2<\infty$ for $z\in K$ where the constants $C_1,C_2$ depend on $K$  but do not depend on $n$ and $z\in K$. Evidently that for such pairs of sequences and functions we have that $|\alpha_n/\beta_n|^{1/n}\to1$ as $n\to\infty$.

Let $\RS_3=\RS_3(\ff_\infty)$ be the three-sheeted Riemann surface associated by Nuttall with the given element $\ff_\infty$, $\ff\in\mathscr Z_{\pm1/2}(E)$, $\ff(\infty)=\prod\limits_{j=1}^m\sqrt{A_jB_j}>0$; see~\cite{Nut84},~\cite{RaSu13},~\cite{Sue18d},~\cite{Sue19}.The zero sheet $\RS_3^{(0)}$ of RS $\RS_3$ is equivalent to Riemann sphere cutted along the segment $E$, $\RS_3^{(0)}\simeq\myh{\CC}\setminus{E}$, and  we have $\partial\RS_3^{(0)}=E^{(0,1)}$, where $\pi_3(E^{(0,1)})=E$. The first sheet $\RS_3^{(1)}$ of $\RS_3$ is equivalent to the Riemann sphere cutted along the compact sets $E$ and $F=\bigsqcup\limits_{j=1}^m[a_j,b_j]$, $\RS_3^{(1)}\simeq\myh{\CC}\setminus(E\sqcup F)$, and have $\partial\RS_3^{(1)}=E^{(0,1)}\sqcup F^{(1,2)}$ where $\pi_3(F^{(1,2)})=F$. The second sheet $\RS_3^{(2)}$ is equivalent to the Riemann sphere cutted along to the compact set $F$, $\RS_3^{(2)}\simeq\myh{\CC}\setminus{F}$, and we have that $\partial\RS_3^{(2)}=F^{(1,2)}$, where $\pi_3(F^{(1,2)})=F$; see Fig.~\ref{fig_1}.

Let prove the inequity~\eqref{7}. We have that
$$
u(z^{(2)})-u(z^{(1)})=2G_F^{\lambda_E}>0,\quad z\notin F.
$$
From~\cite[(16)]{BuSu15} it follows the identity
\begin{equation}
3V^{\lambda_E}(z)+G_F^{\lambda_E}(z)+G_E^{\lambda_F}(z)+3g_E(z,\infty)\equiv w_E=\const,
\quad z\in\myh{\CC},
\label{10}
\end{equation}
where $g_E(z,\infty)$ is Green function for the domain $D$ with the logarithmic singularity at the infinity point.
From~\eqref{10} we obtain that
$$
u(z^{(1)})-u(z^{(0)})=-G_F^{\lambda_E}(z)-3V^{\lambda_E}(z)+w_E
=G_E^{\lambda_F}(z)+3g_E(z,\infty)>0,\quad z\notin{E}.
$$
The relations~\eqref{7} follow directly from the last inequalities.
It the future we will be needed the following equality
\begin{equation}
u(z^{(2)})-u(z^{(0)})=2G_F^{\lambda_E}(z)+G_E^{\lambda_F}(z)
+3g_E(z,\infty)>0, \quad z\notin E\cup F.
\label{11}
\end{equation}

\subsection{}\label{s2s2}
For $\rho\in(1,+\infty)$, let $\Gamma_\rho$ be a level curve of the function $G_F^{\lambda_E}$, i.e.,
\begin{equation}
\Gamma_\rho:=\{z\in\CC:G_F^{\lambda_E}(z)=\log{\rho}\}.
\label{12}
\end{equation}
Since $G_F^{\lambda_E}(z)\equiv0$ when $z\in{F}$ and $[a_k,b_k]\cap[a_j,b_j]=\varnothing$, $k\neq j$, then for some $R>1$ and any $\rho\in(1,R)$
the set $\Gamma_\rho$ consists of exactly $m$ noninterlaced closed components $(\Gamma_\rho)_j$ with the property $\myint(\Gamma_\rho)_j\supset F_j$, $F_j:=[a_j,b_j]$, $j=1,\dots,m$. Set $\Gamma^{(2)}_\rho:=\{z^{(2)}:z\in\Gamma_\rho\}$, $\Gamma^{(1)}_\rho:=\{z^{(1)}: z\in\Gamma_\rho\}$, $\rho\in(1,R)$.
For clarification of the introduced structure see Fig.~\ref{fig_1} below.

Let $V^{(1,2)}\subset\RS_3$ be a neighborhood of the compact set $F^{(1,2)}$ with the property that $\pi_3(\partial V^{(1,2)})=\Gamma_R$ and the element $f_{\infty^{(0)}}$ extends into the domain $\mathfrak D:=\RS_3^{(0)}\sqcup{E^{(0,1)}}\sqcup\RS_3^{(1)}\cup V^{(1,2)}$ as a meromorphic (single-valued) function, $f\in\sM(\mathfrak D)$. From that it directly follows that the error function $R_n(z)$ is also lifted to the point $\infty^{(0)}$ and then extends from this point into the domain $\mathfrak D$ as a meromorphic function $R_n(\zz)$, $\zz\in\mathfrak D$. The function $R_n(\zz)$ has a zero at the point $\zz=\infty^{(0)}$ of order not less then $2n+2$ and a pole at the point $\zz=\infty^{(1)}$ of order not greater then $n$ (we suppose here for simplicity that $f$ is a holomorphic function at the point $\zz=\infty^{(1)}$). Also $f$ can have some other poles at the points in the domain $\mathfrak D$. Let $q_s(z)=z^s+\dots$, where $s\in\NN$ is fixed,  be a polynomial with the following property: the function $q_sf$ is a holomorphic function in the domain $\mathfrak D\setminus\infty^{(1)}$. Let $\mathfrak D_\rho\subset\mathfrak D$ be a domain on RS $\RS_3(f_\infty)$ with the boundary $\partial\mathfrak D_\rho=\Gamma^{(2)}_\rho$, $\rho\in(1,R)$. Up to the end of the paper we will consider only such $\rho\in(1,R)$ that $q_s(z)f(\zz)\neq0$ and $f(z^{(0)})-f(z^{(1)})\neq0$ when $z\in\Gamma_\rho$ (it is easy to see that $f(z^{(0)})-f(z^{(1)})\not\equiv0$). We will refer to such values of $\rho$ as admissible values.

Let $g(\zz):=-u(\zz)$, $\zz\in\RS_3(\ff_\infty)$  (it is the so-called ``$g$-function''; cf.~\cite{IkSu21b},~\cite{Sue21b}). Set (cf.~\cite[(6)]{Sue21d})
\begin{equation}
u_n(\zz):=\log|q_s(z)R_n(\zz)|+(n+1-s)g(\zz),\quad\zz\in\mathfrak D.
\label{14}
\end{equation}
Function $u_n(\zz)$ is a subharmonic function in the domain $\mathfrak D_\rho$ and is a continuous function in a neighborhood of $\Gamma^{(2)}_\rho$, $\rho\in(1,R)$. Therefore by the maximum principle we have that
\begin{equation}
\max_{z\in\Gamma_\rho}|u_n(z^{(2)})|
>\max_{z\in\Gamma_\rho}|u_n(z^{(1)})|,\quad\rho\in(1,R).
\label{15}
\end{equation}
From the definition of $u_n$ it follows the representation
$$
e^{u_n(\zz)}=\bigl|q_s(z)R_n(\zz)e^{(n+1-s)g(\zz)}\bigr|.
$$
Thus from the definition of $g$-function we obtain that
\begin{align}
\max_{\zz\in\Gamma^{(2)}_\rho}e^{u_n(\zz)}
&=\max_{z\in\Gamma_\rho}\Bigl|q_s(z)R_n(z^{(2)})e^{(n+1-s)V^{\lambda_E}(z)-(n+1-s)G_F^{\lambda_E}(z)}\Bigr|
\notag\\
&=
\max_{z\in\Gamma_\rho}\Bigl|q_s(z)R_n(z^{(2)})e^{(n+1-s)V^{\lambda_E}(z)}\Bigr|\frac1{\rho^{n+1-s}},
\label{16}\\
\max_{\zz\in\Gamma^{(1)}_\rho}e^{u_n(\zz)}
&=\max_{z\in\Gamma_\rho}\Bigl|q_s(z)R_n(z^{(1)})e^{(n+1-s)V^{\lambda_E}(z)+(n+1-s)G_F^{\lambda_E}(z)}\Bigr|
\notag\\
&=
\max_{z\in\Gamma_\rho}\Bigl|q_s(z)R_n(z^{(1)})e^{(n+1-s)V^{\lambda_E}(z)}\Bigr|\rho^{n+1-s}.
\label{17}
\end{align}
From~\eqref{15},~\eqref{16} and~\eqref{17} we easy obtain the estimate
\begin{equation}
\max_{z\in\Gamma_\rho}\biggl|q_s(z)R_n(z^{(1)})e^{(n+1-s)V^{\lambda_E}(z)}\biggr|<\frac1{\rho^{2(n+1-s)}}
\max_{z\in\Gamma_\rho}\biggl|q_s(z)R_n(z^{(2)})e^{(n+1-s)V^{\lambda_E}(z)}\biggr|.
\label{18}
\end{equation}
In the similar way by applying the maximum principle  to the function $u_n(\zz)$ we obtain the inequality for $1<\rho<\rho_2<R$:
\begin{equation}
\max_{z\in\Gamma_{\rho_2}}
\biggl|q_s(z)R_n(z^{(2)})e^{(n+1-s)V^{\lambda_E}(z)}\biggr|>\(\frac{\rho_2}{\rho}\)^{n+1-s}
\max_{z\in\Gamma_\rho}
\biggl|q_s(z)R_n(z^{(2)})e^{(n+1-s)V^{\lambda_E}(z)}\biggr|.
\label{19}
\end{equation}

\subsection{}\label{s2s3}
It is easy to check the representation
\begin{align}
R_n(z^{(1)})&=R_n(z^{(2)})+\bigl[Q_{n,1}(z)(f(z^{(1)})-f(z^{(2)}))
+Q_{n,2}(z)\bigl(f^2(z^{(1)})-f^2(z^{(2)})\bigr)\bigr]
\notag\\
&=R_n(z^{(2)})+\bigl(f(z^{(1)})-f(z^{(2)})\bigr)
\bigl[Q_{n,1}(z)+Q_{n,2}(z)\bigl(f(z^{(1)})+f(z^{(2)})\bigr],
\label{20}
\end{align}
where $f(z^{(1)})-f(z^{(2)})\neq0$, $z\in\Gamma_\rho$, for any admissible $\rho\in(1,R)$. From~\eqref{18} and~\eqref{20} it easy follows that for any admissible $\rho$ and as $n\to\infty$ we have that
\begin{equation}
\max_{z\in\Gamma_\rho}
\left|\bigl[Q_{n,1}(z)+Q_{n,2}(z)\bigl(f(z^{(2)})+f(z^{(2)})\bigr)
\bigr]e^{(n+1-s)V^{\lambda_E}(z)}\right|
\asymp
\max_{z\in\Gamma_\rho}
\left|R_n(z^{(2)})e^{nV^{\lambda_E}(z)}\right|.
\label{21}
\end{equation}
From maximum principle applied to the function $u_n(\zz)$ it easy follows the inequality
\begin{equation}
\max_{z\in\Gamma_\rho}\left|R_n(z^{(0)})e^{(n+1-s)V^{\lambda_E}(z)}
\right|\leq\frac1{\rho^{2(n+1-s)}}\max_{z\in\Gamma_\rho}\left|
R_n(z^{(1)})e^{(n+1-s)V^{\lambda_e}(z)}\right|.
\label{22}
\end{equation}
From this inequality and the identity
\begin{equation}
R_n(z^{(0)})=R_n(z^{(1)})+\bigl[f(z^{(0)}) -f(z^{(1)})\bigr]
\[Q_{n,1}(z)+Q_{n,2}(z)\bigl(f(z^{(0)})+f(z^{(1)})\bigr)\]
\label{23}
\end{equation}
we obtain that as $n\to\infty$
\begin{gather}
\max_{z\in\Gamma_\rho}\left|\[Q_{n,1}(z)+Q_{n,2}(z)\bigl(f(z^{(0)})+f(z^{(1)})\bigr)\]e^{(n-s+1)V^{\lambda_E}(z)}\right|
\notag\\
\asymp
\max_{z\in\Gamma_\rho}\left|R_n(z^{(1)})
e^{(n-s+1)V^{\lambda_E}(z)}\right|.
\label{24}
\end{gather}
In the same way as above (i.e., based on the maximum principle applied to the function $u_n(\zz)$) we obtain the following relations as $n\to\infty$
\begin{align}
\max_{z\in\Gamma_\rho}
&\left|
\bigl[Q_{n,1}(z)+Q_{n,2}(z)(f(z^{(0)})+f(z^{(1)})\bigr]
e^{(n+1-s)V^{\lambda_E}(z)}\right|\notag\\
&\qquad\asymp
\max_{z\in\Gamma_\rho}\left|
R_n(z^{(1)})e^{(n+1-s)V^{\lambda_E}(z)}\right|,
\label{25}\\
\max_{z\in\Gamma_\rho}&\left|
\bigl[Q_{n,1}(z)+Q_{n,2}(z)(f(z^{(0)})+f(z^{(2)})\bigr]
e^{(n+1-s)V^{\lambda_E}(z)}\right|\notag\\
&\qquad\asymp
\max_{z\in\Gamma_\rho}\left|
R_n(z^{(2)})e^{(n+1-s)V^{\lambda_E}(z)}\right|.
\label{26}
\end{align}
From~\eqref{18},~\eqref{25},~\eqref{26} and the identity
\begin{align}
Q_{n,1}(z)&+Q_{n,2}(z)(f(z^{(0)})+f(z^{(2)}))
\notag\\
&=Q_{n,1}(z)+Q_{n,2}(z)\bigl(f(z^{(0)})+f(z^{(1)})\bigr)
+Q_{n,2}(z)(f(z^{(2)})-f(z^{(1)}))
\notag
\end{align}
it follows that
\begin{equation}
\max_{z\in\Gamma_\rho}\left|Q_{n,2}(z)e^{(n+1-s)V^{\lambda_E}(z)}
\right|\asymp
\max_{z\in\Gamma_\rho}\left|R_n(z^{(2)})e^{(n+1-s)V^{\lambda_E}(z)}
\right|.
\label{27}
\end{equation}

\subsection{}\label{s2s4}
Since $\lambda_F$ is the balayage of the  measure $\lambda_E\in M_1(E)$ from the domain $\Omega$ to the compact set $F$, we have the identity
\begin{equation}
V^{\lambda_F}(z)=V^{\lambda_E}(z)-G^{\lambda_E}_F(z)+\const,
\quad z\in\myh{\CC},
\label{28}
\end{equation}
with some constant.  Note now that the function $q_s(z)(f(z^{(0)})+f(z^{(1)}))$ is a holomorphic function in the domain $\Omega$.
From this fact and~\eqref{28} it follows that the function
\begin{equation}
v_n(z):=\log\left|q_s(z)
\bigl[Q_{n,1}(z)+Q_{n,2}(z)(f(z^{(0)})+f(z^{(1)}))\bigr]\right|
+(n+s)V^{\lambda_F}(z)
\label{29}
\end{equation}
is a subharmonic function in $\Omega$. Therefore by the maximum principle and taking into account that $\Gamma_\rho$ is a level curve for the Green potential $G_F^{\lambda_E}(z)$, we obtain that for $1<\rho<\rho_2<R$
\begin{gather}
\max_{z\in\Gamma_{\rho_2}}
\left|\bigl[Q_{n,1}(z)+Q_{n,2}(z)(f(z^{(0)})+f(z^{(1)}))\bigr]
e^{nV^{\lambda_E}(z)}\right|
\notag\\
\leq\(\frac{\rho_2}{\rho}\)^n
\max_{z\in\Gamma_{\rho}}
\left|\bigl[Q_{n,1}(z)+Q_{n,2}(z)(f(z^{(0)})+f(z^{(1)}))\bigr]
e^{nV^{\lambda_E}(z)}\right|.
\label{30}
\end{gather}
Since the function
\begin{equation}
\log|Q_{n,2}(z)|+n\bigl(V^{\lambda_E}(z)-G_F^{\lambda_E}(z)\bigr)
\notag
\end{equation}
is a subharmonic function in $\Omega$, then the following inequality holds
\begin{equation}
\max_{z\in\Gamma_{\rho_2}}\left|Q_{n,2}(z)e^{nV^{\lambda_E}(z)}\right|
\leq\( \frac{\rho_2}{\rho}\)^n
\max_{z\in\Gamma_{\rho}} \left|Q_{n,2}(z)e^{nV^{\lambda_E}(z)}\right|.
\label{31}
\end{equation}
Note that the inequality~\eqref{31} is an analog of the Bernshtein--Walsh Theorem but for the level curves of the Green potential $G_F^{\lambda_E}(z)$ instead of the Green function $g_S(z,\infty)$  (see~\cite{Sue21d}).

Combining~\eqref{19},~\eqref{27} and~\eqref{31}, we obtain finally that for any admissible $\rho_2,\rho\in(1,R)$ and as $n\to\infty$ the following asymptotic relations valid
\begin{align}
\max_{z\in\Gamma_{\rho_2}}\left|
Q_{n,2}(z)e^{nV^{\lambda_E}(z)}\right|\asymp
\(\frac{\rho_2}{\rho}\)^n\max_{z\in\Gamma_{\rho}}\left|
Q_{n,2}(z)e^{nV^{\lambda_E}(z)}\right|,
\label{32}\\
\max_{z\in\Gamma_{\rho_2}}\left|
R_n(z^{(2)})e^{nV^{\lambda_E}(z)}\right|\asymp
\(\frac{\rho_2}{\rho}\)^n\max_{z\in\Gamma_{\rho}}\left|
R_n(z^{(2)})e^{nV^{\lambda_E}(z)}\right|.
\label{33}
\end{align}
From~\eqref{32} it follows (cf.~\cite{Sue21d}) that
the sequence $\{n-\mdeg{Q_{n,2}}\}$ is bounded and
\begin{equation}
\frac1n\chi(Q_{n,2})\overset{*}\longrightarrow\lambda_F,
\quad n\to\infty.
\label{34}
\end{equation}
Based on these results it easy to check that
the sequences $\{n-\mdeg{Q_{n,0}}\},\{n-\mdeg{Q_{n,1}}\}$ are  also bounded and
\begin{equation}
\frac1n\chi(Q_{n,j})\overset{*}\longrightarrow\lambda_F,
\quad j=0,1, \quad n\to\infty.
\notag
\end{equation}
Finally from~\eqref{18} and~\eqref{27} we obtain the estimate
\begin{equation}
\varlimsup_{n\to\infty}
\left|R_n(z^{(1)})e^{nV^{\lambda_E}(z)}\right|^{1/n}\leq\frac1{\rho^2},
\quad z\in\Gamma_\rho.
\label{35}
\end{equation}
Now let consider the function $u_n(\zz)$ (see~\eqref{14}) in the open set $V^{(1,2)}(R)$, where $\pi_3(\partial V^{(1,2)}(R))=\Gamma_R$. Function $u_n(\zz)$ is a (pice-wise) subharmonic function in $V^{(1,2)}(R)$. Ultimately from~\eqref{33},~\eqref{34},~\eqref{35} and by the three constants theorem applied to $u_n$ and for any component of $V^{(1,2)}(R)$ we obtain that as $n\to\infty$
\begin{equation}
\max_{z\in\Gamma_\rho}\left|R_n(z^{(1)})e^{nV^{\lambda_E}(z)}\right|\asymp\frac1{\rho^2}.
\label{36}
\end{equation}
From this relation,~\eqref{25} and~\eqref{34} it finally follows that
\begin{equation}
\biggl|\frac{Q_{n,1}(z)}{Q_{n,2}(z)}+\bigl(f(z^{(0)})+f(z^{(1)})\bigr)
\biggr|^{1/n}\overset{\mcap}\longrightarrow e^{-2G_F^{\lambda_E}(z)}<1,\quad z\in\Omega\setminus{E}.
\label{37}
\end{equation}
Theorem~\ref{the1} is proven.

\section{Concluding Remarks and Some Conjectures}\label{s3}

\subsection{}\label{s3s1}
Let $f\in\mathbb C(z,\ff)$, $\ff_\infty$ be the element of the function $\ff$ specified above and $f_\infty\in\HH(\infty)$ be the corresponding element of $f$. Let $n\in\NN$ and $P_{2n},P_{2n,1},P_{2n,2}$, $\mdeg{P_{2n}},\mdeg{P_{2n,1}},\mdeg{P_{2n,2}}\leq{2n}$, be the corresponding type II Hermite--Pad\'e polynomials for the pair of functions $f,f^2$ and multiindex $(2n,2n)$, i.e.
\begin{equation}
\begin{aligned}
R_{n,1}(z):&=(P_{2n}f_\infty-P_{2n,1})(z)=O\(\frac1{z^{n+1}}\),
\quad z\to\infty,\\
R_{n,2}(z):&=(P_{2n}f_\infty^2-P_{2n,2})(z)=O\(\frac1{z^{n+1}}\),
\quad z\to\infty.
\end{aligned}
\label{3.1}
\end{equation}
It would be natural to conjecture that the following statement is valid.

\begin{conjecture}\label{conj1}
If $f\in\CC(z,\ff)$ and $f_\infty\in\HH(\infty)$ corresponds to the $w_\infty$ specified  above, then we have as $n\to\infty$
\begin{align}
\frac1n\chi(P_{2n}),&\frac1n\chi(P_{2n,1}),\frac1n\chi(P_{2n,2})
\overset{*}\longrightarrow\lambda_E,\label{3.2}\\
\frac{P_{2n,1}}{P_{2n}}(z)&\overset\mcap\longrightarrow f(z^{(0)}),
\quad \text{inside}\quad D\label{3.3}\\
\frac{P_{2n,2}}{P_{2n}}(z)&\overset\mcap\longrightarrow f^2(z^{(0)}),
\quad \text{inside}\quad D\label{3.4}.
\end{align}
\end{conjecture}

Note that since $f\in\mathbb C(z,\ff)$ in general is a complex-valued function on the real line, then the powerful methods developed in the papers~\cite{GoRa81},~\cite{GoRaSo97} and~\cite{ApLy10} are not applicable to prove the relations~\eqref{3.2}--\eqref{3.4} (in this connection see also~\cite{MaRaSu16},~\cite{ApLoMa17},~\cite{Rak18}, \cite{Sor20},~\cite{ApLy21}). If~\eqref{3.3} would be proven, then from that and~\eqref{37} it would be followed that on the base of type I and type II HP-polynomials a multi-valued analytic function $f\in\mathbb C(z,\ff)$ is constructively recovered on the two sheets of the RS $\RS_3(\ff_\infty)$. Here we use the term ``constructive recovering'' in the sense of the paper of P.~Henrici~\cite[Sec.~2]{Hen66} (on that subject see~\cite{BoBoSh06},~\cite{ApBuMaSu11},~\cite{Sue15},~\cite{KoPaSuCh17},~\cite{Sue18d},~\cite{Sor20},~\cite{Bus20},~\cite{IkSu21},~\cite{Kom21} and the bibliography therein).

Note that so far all the consideration in this paper were connected with the element $\ff_\infty=\ff_{\infty^{(0)}}$. Therefore the next question arises in a natural way when we consider as a start point the element $\ff_{\infty^{(1)}}$ but not the element $\ff_{\infty^{(0)}}$ as it was done in this paper. What will be there a three-sheeted RS  $\RS_3(\ff_{\infty^{(1)}})$ associated by Nuttall with that element $\ff_{\infty^{(1)}}$ and what can we say about the asymptotics of HP-polynomials associated with the chosen element $\ff_{\infty^{(1)}}$ and with the corresponding element $f_{\infty^{(1)}}$?

\subsection{}\label{s3s2}
Traditionally the problem of asymptotics of Pad\'e and Hermite--Pad\'e polynomials consists of two components. Namely, geometrical component and analytical one. In the current paper we suppose the geometrical component to be trivial, i.e. both plates of the Nuttall condenser are located on the real line, $E,F\subset\RR$. But since $f\in\CC(z,\ff)$ is a complex-valued function on the real line, then the analytic component is not trivial. Thus the next step is very natural.

\begin{figure}
\centerline{
\boxed{
\includegraphics[width=\textwidth]{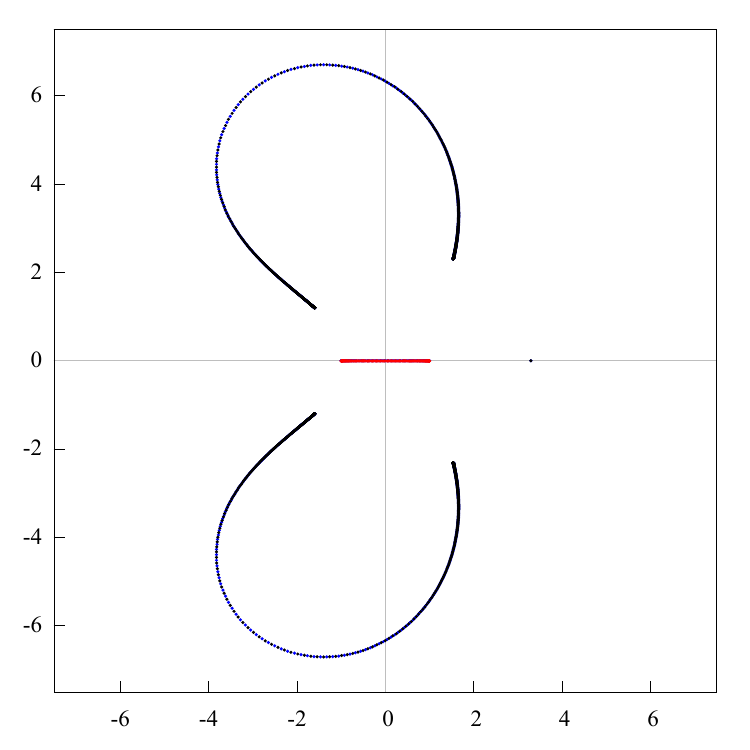}
}}
\caption{
Here the results of numerical simulations are represented. We consider the case when in~\eqref{1} $m=2$ and $A_1=\myo{B}_1,A_2=\myo{B}_2\not\in\RR$. The compact set $E$ is still  the segment $[-1,1]$ but  now the compact set $F$ is symmetric with respect to real line and consists of two arcs. The red points are the zeros of type I
HP-polynomial $P_{100}$. They simulate the compact set $E$. The black points are the zeros of type II HP-polynomial $Q_{1200,2}$. They simulate the compact set  $F$. }
\label{fig_2}
\end{figure}

Let suppose that in the representation~\eqref{1} a more general situation is admissible, i.e. for some $j\in\{1,\dots,m\}$ we still have that $A_j<B_j$ but some other $j\in\{1,\dots,m\}$ we have that $A_j=\myo{B}_j\notin\RR$. This case requires to consider the general definition of the Nuttall condenser $(E,F)$ as it was done by E.~A.~Rakhmanov and the author in~\cite{RaSu13} in 2013 (see also the Fig.~\ref{fig_2}). Namely, the situation when  the plate $E\subset\RR$ but it is not the case for the plate $F$. Instead it is supposed that $F$ is symmetric with respect to the real line, i.e. $z\in F$ $\Leftrightarrow$ $\myo{z}\in F$. Of course in any case it is supposed that $E\cap F=\varnothing$.

\begin{conjecture}\label{conj2}
If $f\in\CC(z,\ff)$ and $f_\infty\in\HH(\infty)$ then under the assumptions of this subsection~\ref{s3s2} we have the relations~\eqref{8}--\eqref{9} and~\eqref{3.2}--\eqref{3.4} are valid  as $n\to\infty$.
\end{conjecture}

\vskip5mm
\newpage

\vskip5mm

\def\by#1;{#1\unskip,}
\def\paper#1;{``#1\unskip''\unskip,}
\def\paperinfo#1;{#1\unskip.}
\def\eprint#1;{``#1\unskip''\unskip,}
\def\eprintinfo#1;{#1\unskip,}
\def\book#1;{``#1\unskip''\unskip,}
\def\inbook#1;{``#1\unskip''\unskip,}
\def\bookinfo#1;{#1\unskip,}
\def\jour#1;{#1\unskip,}
\def\issue#1;{#1\unskip,}
\def\yr#1;{#1\unskip,}
\def\pages#1.{#1\unskip.}
\def\vol#1;{\textbf{#1}\unskip,}
\def\finalinfo#1;{#1\unskip.}
\def\publ#1;{#1\unskip,}
\def\publadrr#1;{#1\unskip,}
\def\publaddr#1;{#1\unskip,}
\def\procinfo#1;{#1\unskip,}
\def\serial#1;{#1\unskip,}
\def\ed#1;{ed #1\unskip,}
\def\eds#1;{ed #1\unskip,}

\end{document}